\documentclass{jnmp}

\usepackage{amsmath}

\newtheorem*{theorem}{Theorem}

\begin{document}

\renewcommand{\evenhead}{M Schwarz Jr}
\renewcommand{\oddhead}{Nonlinear Schr\"odinger, Infinite Dimensional Tori and Neighboring Tori}

\thispagestyle{empty}

\FirstPageHead{10}{1}{2003}{\pageref{schwarz-firstpage}--\pageref{schwarz-lastpage}}{Article}

\copyrightnote{2003}{M Schwarz Jr}

\Name{Nonlinear Schr\"odinger, Infinite Dimensional\\
 Tori and Neighboring Tori}
\label{schwarz-firstpage}

\Author{M SCHWARZ Jr}

\Address{Mathematics Department, Northeastern University,
Boston Massachusetts 02115, USA\\
}

\Date{Received March 23, 2002; Accepted July 17, 2002}

\begin{abstract}
\noindent
In this work, we explain in what sense the generic level set of the constants of
motion for the periodic nonlinear Schr\"odinger 
equation is an infinite dimensional torus on which each generalized nonlinear
Schr\"odinger flow is reduced to straight line almost periodic motion, and describe how
neighboring generic infinite dimensional tori are connected.
\end{abstract}

\section{Introduction}
We consider the Hamiltonian equation 
\begin{equation}
iu_t+u_{xx}-|u|^2u=0
\end{equation}
 of the periodic nonlinear Schr\"odinger equation, where
$u(x,t)$ is a complex valued function in the class of smooth period one
functions. In this work, we explain in what sense the generic level set of the constants of
motion for the periodic nonlinear Schr\"odinger 
equation is an infinite dimensional torus, why the
solution of the Hamiltonian equation is almost periodic in time, and describe how
neighboring generic infinite dimensional tori are connected. Bourgain~[1] has solved the initial
value problem for the periodic nonlinear Schr\"odinger equation. Ma and Ablowitz~[2] 
have reduced the
periodic nonlinear Schr\"odinger equation to an inverse spectral problem for periodic potentials.
They provide explicit formulas for the special class of $N$-soliton solutions of the periodic
nonlinear Schr\"odinger equation and found an infinite sequence of functionals that are in
involution and constant along solutions of~(1). For the nonlinear Schr\"odinger 
equation, Batig et al~[3] and Schmidt~[4] used 
the method of inverse spectral theory and integrated the equation in
the class of analytic~[4] and smooth periodic functions~[3]. They identified the generic invariant
set of the constants of motion with an infinite dimensional tori. Their study~[2, 3] did not
describe how neighboring tori are connected. 

The nonlinear Schr\"odinger equation is an example of the Hamiltonian equation  
\begin{equation}
\frac{\partial u}{\partial t}=K(u),
\end{equation}
 where $K(u)$ is a nonlinear
operator and $u$ a complex valued function in the class of smooth periodic
functions. Let $F_m(u)$ denote functionals that are in involution and constant
along solutions of (2). In~[5], we give a proof of an infinite dimensional
version of Liouville's theorem and explain in what sense the generic level set of the
functionals $F_m(u)$ is an infinite dimensional torus on which the solution of~(2)
reduces to straight line motion that is almost periodic in time. Furthermore, we explain in what
sense neighboring generic tori and solutions of (2) are connected. The approach in~[5] is related
to Lax's~[6] study of finite-dimensional level sets of completely integrable partial differential
equations and is independent of the method of inverse spectral theory and the viewpoint of
algebraic curves. An application of the theorem in~[5] to the nonlinear Schr\"odinger equation
yields a different proof of the result of Batig~[3] and Schmidt~[4]. In addition, the present work
describes how neighboring generic tori and the solutions of (2) are connected.

In the classical case
\[
\frac{dv}{dt}=K(v),\qquad v\in {\mathbb R}^{2N}
\]
a theorem of Liouville [7] states that the system is completely integrable. If the involutive
constant functions $F_m(v)$,
$m=1,2,\ldots,N$ are independent in the sense that their gra\-dients are
linearly independent and if the $N$ dimensional level set satisfying
$F_m(v)=F_m(v_0)$, $m=1,2,\ldots,N$ is compact; in fact,
\begin{enumerate}
\itemsep0mm
\vspace{-2mm}
\item[(a)] the level set is an $N$ dimensional torus on which the
flow is quasiperiodic and  
\item[(b)] neighboring Louville tori are diffeomorphic to one another.\vspace{-2mm}
\end{enumerate}

The proof of the classical Liouville theorem is based on the inverse function theorem. It verifies
that the composition of the commuting flows associated with $F_m(v)$,
$m=1,\ldots,N$ identifies a
neighborhood of
${\mathbb R}^N$ with a neighborhood of the level set. 
The basic periods of this map are used to identify a
connected component of the level set with an $N$ dimensional 
torus on which the Hamiltonian flow
associated with each $F_m(v)$ is reduced to straight line quasiperiodic motion.

In [5], we gave a proof of an infinite dimensional version of Liouville's theorem.
We were unable
to use the inverse function theorem. We introduced instead a local open
mapping theorem for certain
types of nondifferentiable maps and established that the composition of the commuting flows
associated with $F_m(u)$,
$m\ge 1$ defines a continuous open map from the Hilbert space $l_2$ of
square summable sequences onto a connected component of a
generic compact level set. This map is
not locally diffeomorphic because $l_2$ in not locally compact. The periods of this mapping
are contained in any neighborhood of $l_2$. A complete set of 
basic periods was used to identify a
connected component of the level set with an infinite dimensional
torus on which the Hamiltonian
flow associated with each $F_m(u)$ is reduced to straight line almost periodic motion. 
Furthermore
we established that the complete set of basic periods 
that characterized a generic level set may
be continuously extended to a complete set of basic periods that describe a neighboring
generic level set. We established a sense in which neighboring generic 
level sets are homeomorphic
to the standard infinite dimensional torus, and determined a sense in 
which these neighboring level
sets are connected. This present study of the periodic nonlinear Schr\"odinger equation is an
illustration of the result~[5].

\section{Theorem}

 Let $W_n$ $(n\ge 0)$ denote the usual Sobolev space of
functions on $[0,1]$, of period one, having derivatives of all orders up to $n$
with norm 
\[
\|w\|_n^2=\sum_{j \le n}
\int_0^1   |D^jw(x)|^2 dx.
\]
  The norm in the space $L_2$ is denoted
by $\|w\|$. For $w\in W_n$ and integers $j$, $k$, and $p$ with
$p \ge 2$, it is known that   
\[
\sqrt[p]{\int_0^1|D^jw(x)|^pdx} \le
2^{p-2/2p} \|D^kw\|^a \|w\|^{1-a},
\]
  where $a=\left(j+\frac 12 - \frac 1p\right)/k$ and $1 \le j<k \le n$. 
We denote by $C_1^n$ the space of functions
of period one having continuous derivatives of order less than or equal $n$. 
The value of $n$ does
not enter into the proof of the result [5]. The value of $n$ specifies the class of solutions for
the nonlinear Schr\"odinger equation or generalized equations.
The subscript of $W_n$ is generally suppressed.

  The Hamiltonian formulation of (2) is due to Gardiner [8] and Lax [6]. Let
$F(u)$ denote a functional whose argument is a smooth function of period one and let 
$(\cdot ,\cdot )$
denote the scalar product in $L_2$.  Then  
\[
\lim_{\epsilon \rightarrow0} \epsilon^{-1}(F(u+\epsilon v)-F(u))=(G_F(u),v)
\]
 for appropriate $u$ and $v$
defines $G_F(u)$, the gradient of $F$ at $u$. Define the Poisson bracket of
$F(u)$ with $H(u)$ by 
\[
\{ F(u),H(u) \}=(G_F(u),JG_H(u)),
\]
  where $J$
is an antisymmetric operator independent of $u$.
If ${\rm K}(u)=JG_F(u)$, then the equation (2) is
said to be Hamiltonian. We denote by $S_F(t)u$ the nonlinear operator determining the
solution of (2) on the basis of its initial values at  $t=0$: $u(t)=S_{F}(t)u_0$.
If $\{ F,H\}$ =0 for all $u$, then the solutions of (2) and of $u_t=JG_H(u)$
commute: $S_H(t)S_F(t')=S_F(t')S_H(t)$ for all $t$ and $t'$. 

As for the nonlinear Schr\"odinger equation, Ma and Ablowitz have constructed
explicitly an infinite sequence of functionals $I_m(u)$ that are constant along
the flow (1).  The first three are 
\[
\int_0^1u\overline{u}\,dx,\qquad
\frac i2 \int_0^1(\overline{u} u_x-u\overline{u}_x)\,dx,\qquad
\int_0^1\left(|u_x|^2+\|u|^4\right).
\]
 Let $G_{I_m}(u)$ denote the gradient of
$I_m$ with respect to $\overline{u}$ at $u$, and let
\[
\langle u,v\rangle =\int_0^1(u\overline{v}+v\overline{u})\,dx
\]
 denote the product in $L_2$. Then the Poisson bracket 
\[
\{I_m(u),I_n(u)\}=\langle G_{I_m}(u),JG_{I_n}(u)\rangle
\]
where the symplectic structure is introduced through $J=i$. It is known that
\[
\{I_m(u),I_n(u)\}=0
\]
 for all $m$ and $n$ and smooth periodic
functions $u$. Therefore $I_m(u)$ are constant along solutions of
\[
\frac{\partial u}{\partial t}=K_{I_m}(u)=JG_{I_m}(u),\qquad m\ge 1
\]
 the generalized nonlinear Schr\"odinger equation, where $m=3$ is equation (1). The
non\-li\-near Schr\"odinger equation (1) is of the form 
\[
u_t=[L,A],
\]
 where the operator L and A depend on $u$. The nonlinear Schr\"odinger flow preserves the
spectrum of $L$ determined by 
\begin{equation}
Lf=\begin{pmatrix}-D&u\cr -\overline{u}&D\end{pmatrix}
\begin{pmatrix}f_1\cr
f_2\end{pmatrix}=i\lambda\begin{pmatrix}f_1\cr f_2\end{pmatrix}
\end{equation}
 in the class of functions $f[0,1)\rightarrow C^2$ with $f(x+1)=mf(x)$ for $0\le x <1$. The
periodic and antiperiodic spectra $[m=\pm 1]$ will be of special interest. Use
of (3) and a direct calculation shows that the spectrum of $L$ is real. Ma and
Ablowitz determined that the periodic and antiperiodic spectrum of $L$ is
comprised of simple and double eigenvalues $\lambda_m$ with eigenfunctions
$f=(f_{1,m},f_{2,m})^T$, and that the functionals $\lambda_m(u)$ are in involution. This study
concerns the general situation in which the spectrum of $L$
 is simple. The exceptional case of mixed
simple and double spectra offers no additional technical difficulties. Let $M$ denote the
portion of the space $W$ of smooth periodic functions for which the spectrum of
$L$ is simple. For
$u_0\in M$, we consider the level set
\[
M_{\lambda(u_0)}=\{u\, |\, \lambda_m(u)=\lambda_m(u_0),\; m\ge 1\}
\]
 in $W_n$. We prove that
$M_{\lambda}$ is generated by the sequence of the generalized 
nonlinear Schr\"odinger flows and that
the generic level set is identified with an infinite dimensional torus on which each generalized
nonlinear Schr\"odinger flow is reduced to straight line motion that is almost periodic in time.
Furthermore, we make precise the sense in which neighboring generic level sets are connected.

 To identify
$M_{\lambda{u_0}}$ with the standard infinite-dimensional torus $T^{\infty}=[0,1)^{\infty}$ we
first state the result in [5]. Consider the Hamiltonian equation (2), the sequence $F_m(u)$ that
are in involution and constant along solutions of (2), and the level set
\[
M_{u_0}=\{u\,|\,F_m(u)=F_m(u_0), \; m\ge 1\}.
\]
Let $u$ be an element of $M_{u_0}$ and view the latter as a subset of $L_2$.
Define $G_{F_m}(u)$ to be the gradient of $F_m(u)$ at $u$. $G_{F_m}(u)$ is a
vector that is normal to $M_{u_0}$ at $u$. Let $N_u$ be the closure in $L_2$ of
the span of $G_{F_m}(u)$ and assume that $G_{F_m}(u)$ is a basis of $N_u$; by
which we mean a) each element $G_u$ in $N_u$ is uniquely expressible as
$G_u=tG(u)=\sum\limits_{m=1}^{\infty}t_mG_{F_m}(u)$ for $t$ in the Hilbert space
$l_2$ and b) $G_u$ admits the estimate  
\[
c_1(u)|t|_{l_2}\le \|tG(u)\|\le
c_2(u)|t|_{l_2},
\]
 where $c_1$ and $c_2$ depend continuously on $u$ in $M$.
$N_u$ is the normal space of $M_{u_0}$ at $u$ and, by our assumptions, no single
gradient $G_{F_m}(u)$ lies in the closure in $L_2$ of the other gradients
$G_{F_n}(u)$. The Poisson bracket of $F_m(u)$ and $F_n(u)$ vanishes for all
$m$ and $n$ and for $u$ in the class of smooth period one functions. The
functionals $F_m(u)$ generate commuting flows
\begin{equation}
\frac{\partial u}{\partial t}=K_{F_m}(u)=JG_{F_m}(u),\qquad   m \ge 1
\end{equation}
 on $M_{u_0}$
and $F_1(u),\ldots,F_m(u),\ldots$ are constants of these motions. $K_{F_m}(u)$ is
tangent to $M_{u_0}$ at $u$. Denote by $T_u$ the closure in $L_2$ of the span of
$K_{F_m}(u)$. Suppose that $K_{F_m}(u)$ is a~basis of $T_u$; each element $K_u$
in $T_u$ is uniquely expressible as
$K_u=\sum\limits_{m=1}^{\infty}t_mK_{F_m}(u)=tK(u)$ for $t$ in $l_2$, and
\begin{equation}
c_1(u)|t|_{l_2}\le \|tK(u)\|\le c_2(u)|t|_{l_2},
\end{equation}
where $c_1$ and $c_2$ depend continuously on $u$ in $M$. Assume that $T_u$ equals the
orthogonal complement of $N_u$. $T_u$ represents the tangent space and every
direction of $L_2$ has been accounted for. Let $S_{F_m}(t_m)u_0$ denote the
nonlinear operator uniquely determining the solution of (4) on the basis of its
initial values at $t=0$: $u(t)=S_{F_m}(t_m)u_0$. For $t$ in $l_2$ we show that 
\[
S(t)u=\lim_{N\rightarrow\infty} \prod_{m=1}^N S_{F_m}(t_m)u_0
\]
  in $W_n$
where $S(t+t')u=S(t)S(t')u$ for $t$, $t'$  in $l_2$, and for $t\in l_2$, 
$S(t)u\in W_m$ is continuous in $t$ uniformly in $u$ on $M_{u_0}$.
Denote by $dG_F(u)$ the second derivative of $F$
defined by 
\[
\lim_{\epsilon \rightarrow 0} \epsilon^{-1}(G_F(u+\epsilon v)-G_F(u))=dG_F(u)v.
\]
Let $v(\tau)$, $\tau\ge 0$ be a curve in $M_{u_0}$ that
satisfies $\frac{dv(\tau)}{d\tau}=K_{v(\tau)}$ with $v(0)=v_1$ and let
$dG_{v(\tau)}K_{v(\tau)}=\frac{dG_{v(\tau)}}{d\tau}$ admit the
estimate
\begin{equation}
\frac{(dG_uK_v, K'_v)}{\|K_v\| \|K'_v\|} \le c\|G_v\|,
\end{equation}
 where $K'_v \in T_v$, $c$ is independent of $v(0)\in
M_{u_0}$, and $v=v(\tau)$ for small $\tau$. Then $S(t)u_0$ is an open map of
$l_2$ onto $M_{u_0}$ in $W_n$. Let $L_{u_0}$ denote the set of $t$ in $l_2$ for
which $S(t)u=u$ for all $u$ in $M_{u_0}$. $S$ is a homeomorphism of $l_2
/L_{u_0}$ onto $M_{u_0}$ in $W_n$. $l_2/L_{u_0}$ is compact and may be
identified as in [5] with the standard infinite-dimensional torus $T^{\infty}$:
in more detail, there exist $\omega_m$, $m\ge 1$ from $L_{u_0}$ so that each $t$
of $l_2/ L_{u_0}$ is uniquely represented by
$t=\sum\limits^{\infty}_{m=1}\tau_m\omega_m$, where $0\le\tau_m< 1 $ for all $m$.
$M_{u_0}$ is an infinite-dimensional torus and the solution $S_{F_m}(t_m)u$ is
almost periodic on $l_2/L_{u_0}$, uniformly with respect to initial values $u\in
M_{u_0}$. The motion $S_{F_m}(t_m)u$ 
of each Hamiltonian equation related to $F_m$ is identified
with straight line motion on $l_2/L$: in detail, for $e_m$ in the $m$-th
coordinate direction in $l_2/L$ and $\hat\omega_m\in l_2$ with
$\omega_n\cdot\hat\omega_m=\delta_{n,m}$, then $S_{F_m}(t_m)u$ is identified with straight
line motion in the direction $\sum\limits_{m=1}^\infty(e_m\cdot\hat\omega_m)\omega_m$ on
$l_2/L$.

  For the generalized nonlinear Schr\"odinger flow,
$K_{I_m}(u)$ is an element of the closure in~$L_2$ of the 
span of $K_{\lambda_m}(u)$ and the flow
of each generalized nonlinear Schr\"odinger equation is identified with straight line motion
that is almost periodic in time on $l_2/L$.

 We next identify as in [5], neighboring generic tori and then state a sense in which they are
related. Let $F_m(u)$ denote a sequence of analytic functions of $u\in W$
and suppose $\sum\limits_{m=1}^{\infty}F_m(u)^2$ is bounded uniformly in u on bounded
sets in $W$. Let $F(u)=(F_1,\ldots,F_m,\ldots)$ and $l=F(M).$ Let $v_0\in M$ and
$f_0=F(v_0)\in l$. Write $u_{f_0}$ for $v_0$. For $f$ in a small neighborhood of
$f_0$, there exists $u_f$ in $M$ and $f=F(u_f)$ and 
\[
| f-f_0|_{l_2}\ge c\|u_f-u_{f_0}\|,
\]
 where $c$ is locally independent of $u\in M$ and $f\in l$.
The curve that joins $u_{f_0}$ to $u_f$ depends uniquely on $u_{f_0}$. The curve
is relatively short in the sense that the length of the curve in $W$ joining
$u_{f_0}$ with $u_f$ is bounded by a fixed multiple of $|f-f_0|_{l_2}$.  The
torus $M_f=F^{-1}(f)$ is homeomorphic to the standard torus $T^{\infty}$.
$M_{u_{f_0}}$ is characterized by basic generators $\omega_m(f_0)$ that are the
periods of $S(t)u_{f_0}$, and for $f$ in a small enough neighborhood of $f_0$ in~$l$, 
the $\omega_m(f_0)$ may be continuously extended to the basic generators
$\omega_m(u_f)$ that describe $M_{u_f}$. $M_{u_f}=F^{-1}(f)$ is identified with
the set $T_f$ of convergent sums $\sum\limits_{m=1}^{\infty}\tau_m\omega_m(f)$, 
$0\le \tau_m < 1$, which converge in $l_2$ uniformly in $\tau_m$ and $f$.
Furthermore, there exists a curve that is continuous in $l_2$ that connects
$\sum\limits_{m=1}^{\infty}\tau_m\omega_m(f_0)$ with
$\sum\limits_{m=1}^{\infty}\tau_m\omega_m(f)$ and is relatively short in the sense that
the length of the curve in $l_2$ is less than $\frac{1}{2}
\left| \sum\limits_{m=1}^{\infty}\tau_m\omega_m(f_0)\right|_{l_2}$. $T_f$ is homeomorphic to
$T_{f_0}$ and $T_f$ is uniformly close to $T_{f_0}$. 
$M_{u_f}$ is homeomorphic to $M_{u_{f_0}}$ and
to the standard infinite-dimensional torus $T^{\infty}$. Furthermore, $M_{u_f}$ and
$M_{u_{f_0}}$ are connected by a relatively short continuous curve in $W$ that
is contained in $M$ except for a countable number of elements. This leads to the
result of this work on the nonlinear Schr\"odinger equation.

\begin{theorem} Let $u\in M$ and $F_m(u)=\lambda_m(u)-\lambda_m(0)$. $F_m(u)$ is a
sequence of analytic functions of $u\in W$ that are in involution and the level
set  $M_u$ is bounded. The sequence $G_m(u)$ and $K_m(u)$ is a
basis for $N_u$ and $T_u$ respectively with $N_u\oplus T_u=L_2$ and $dG$ admits
the estimate $(dG_uK_u,K_uÕ)/\|K_u\|\|K_uÕ\|\le c\|G_u\|$, where $c$ is
independent of $u$. For $u_0\in M$, $S(t)u_0$ is a homeomorphism of
$l_2/L_{u_0}$ onto $M_{u_0}$ in $W$. $l_2/L_{u_0}$ is compact and identified
with an infinite dimensional torus: there exists a sequence $\omega_m$ from
$L_{u_0}$ for which each element of $l_2/L_{u_{f_0}}$ is uniquely represented by
$\sum\limits_{m=1}^{\infty}\tau_m\omega_m$,
$0\le \tau_m<1$ for all $m$. The flow of each generalized
nonlinear Schr\"odinger equation
 is identified with straight motion that is almost periodic in time
on $l_2/L_{u_{f_0}}$. For $u$ in $M$, $F_m(u)$ is square summable uniformly in 
$u$ on bounded sets in
$W$. For directions $v$ transverse to $M_u$ at $u$, $dG_m(u)$ admits the estimate
$\|dG_m(u)v\|\le c_m\|v\|$,
 where $c_m$ is square summable independently of $u$ and $v$. 
For~$f$ in a small neighborhood of
$f_0=F(u_{f_0})$ in $l$ there exists $u_f\in M$ satisfying $F(u_f)=f$ that
admits the estimate $|f-f_0|_{l_2}\ge c\|u_f-u_{f_0}\|$, where c is
locally independent of $u$ in $M$ and $f\in l$. $M_{u_f}$ is homeomorphic to
$M_{f_0}$ and to the standard infinite-dimensional torus. Furthermore $M_{u_f}$
and $M_{u_{f_0}}$ are connected by a relatively short continuous curve in $W$
that is contained in $M$ except for a countable number of elements.
\end{theorem}

This completes the statement of the Theorem.

A modification of the proof of Lemma 1 in [5] 
establishes that $M_{u_f}$ and $M_{u_{f_0}}$ are
connected by a relatively short smooth curve in $W$ that is contained in $M$ except for a
countable number of elements. This verifies that $M_{u_f}$ is diffeomorphic to $M_{f_0}$ and to
the standard infinite-dimensional torus. The exceptional case 
of mixed simple and double spectra
can be solved by a slight modification of the approach taken in [5]. The application of [5] to the
Hamiltonian flow 
\[
iu_t+u_{xx}+|u|^2u=0
\]
 offers no additional technical problem. In this case
the torus is of lower dimension.  

\section{Proof}

In this section we prove the Theorem. We establish
first the properties of $M_{u_0}$ and obtain an apriori estimate of $dG$. For
$u$ in $M$, we verify that $G_m(u)$ is a basis for $N_u$ and $K_m(u)$ is a basis
for $T_u$ and that $N_u\oplus T_u=L_2$.

\medskip

\underline{\it Item 1.} The functionals $\lambda_m(u)$ are in involution
and $M_{u_0}$ is bounded in $W$. For $u_0\in M$, consider 
\[
M_{u_0}=\{u\, |\, F_m(u)=\lambda_m(u)-\lambda(0)=F_m(u_0), \; m\ge 1\}
\]
 and use the result in [2] or the periodic version of the result of Zakharov and
Shabat [9] to show that $I_m(u)$ are directly related to $\lambda_m$ and that
$I_m(u)$ is constant on $M_{u_0}$.  The functional 
\[
I_1=\int_0^1|u|^2dx
\]
 gives 
\[
\|u\|\le c.
\]
 Rearrange the functional 
\[
I_3=\int_0^1\left(|u_x|^2+|u|^4\right)dx
\]
 and estimate to find that
\[
\int_0^1|u_x|^2 dx \le c+c\|u_x\|\|u\|^3 \le 
c+c\|u_x\| \le c+\frac{1}{2}\|u_x\|^2,
\]
 where we have applied the
general inequality $|u|_{\infty}\le \sqrt{2}\|u_x\|^{1/2}\|u\|^{1/2}$ 
and previous bounds. This leads to the estimate
$\|u_x\|\le c$.  The integral 
\[
I_5=\int_0^1|u_{xx}\|^2+\frac{1}{2} |u|^6-\frac{1}{2}
\left(\frac{d}{dx}|u|^2\right)^2-3|u_x|^2|u|^2dx
\]
 is estimated as follows: 
\[
\|u_{xx}\|^2\le c\|u_x\|^{1/3}\|u\|^{2/3}+c|u|_\infty^2\|u_x\|^2\le c,
\]
 where we
have used the general estimate $|u|_6\le 2^{1/3}\|u_x\|^{1/3}
\|u\|^{2/3}$ and previous bounds to find 
\[
\| u_{xx}\|\le c.
\]
 For
$n\ge 4$, the functionals $I_n$  have weight $2n$ where the weight is a sum of
the weights of its factors and the weight of $D^ru$ is $1+r$. Use previous
estimates to obtain 
\[
\|u\|_n\le c
\]
 for any $n$. $M_{u_0}$ is bounded in
$W_n$. 

\medskip

\underline{\it Item 2.} For simple eigenvalue $\lambda$, the gradient of
$\lambda_m$ with respect to $\overline{u}$ equals
\[
\frac{d\lambda_m}{d\bar u}=if_1\overline{f}_2=if_1^2.
\]

 Begin with the equation (3) for $f_1$ and $f_2$. Let
\[
(u,v)=\int_0^1uv\,dx.
\]
 Let $u^{\epsilon}=\overline{u}+\epsilon v$ and compute the
derivatives 
\[
\frac{d}{d\epsilon}\lambda\qquad {\rm and}\qquad \overline{\frac{d}{d\epsilon}\lambda}.
\]
 Begin with the equation for $f_2$ in (3)
and compute the derivative with respect to $\epsilon$. Multiply the resulting
equation by $\overline{f}_2$ and integrate with respect to $x$ from zero to one.
Integrate by parts and use the equations once again and substitute
 $(\dot f_{2x},\overline{f}_2)$ and $(\dot{\overline{u}}_{2x},f_2)$ into the previous
expression, and find 
\begin{gather*}
i\dot{\overline{\lambda}} (\overline{f}_2,f_2)=(\overline{\dot f}_2,\overline{u} f_1+i\lambda
f_2)+(\dot u, f_2\overline{f}_1)+(u,\overline{\dot f}_1 f_2)-i
\overline{\lambda}(\overline{\dot f}_2,f_2),\\
i\dot{\lambda}(f_2,\overline{f}_2)=-(\dot f_2,u\overline{f}_1-i\lambda
\overline{f}_2)-(\dot{\overline{u}},f_1\overline{f}_2)-
(\overline{u},\overline{f}_2\dot{f}_1)-i\lambda(\dot{f}_2,\overline{f}_2).
\end{gather*}
The properties $\overline{\lambda} =\lambda$,
$f_1=\overline{f}_2$, and $f_2=\overline{f}_1$
follow directly from (3). Substitute these identities into
the previous equations with $\|f_1\|=\|f_2\|=1$ and obtain 
\[
\dot{\overline{\lambda}}
+ \dot{\lambda}=(\dot{\overline{u}},if_1\overline{f}_2 )+
(\dot{u},\overline{i\overline{f}_2f_1}).
\]
 Use the inner
product $\langle u,v\rangle$ and find that 
\[
\frac{d\lambda}{d\overline{u}}=if_1\overline{f}_2 =if_1^2.
\]

\underline{\it Item 3.} For $u$ in $M$ and $v$ transverse to $M_u$, 
let $u(\tau)=\tau v + (1-\tau)u$. Then 
\[
| f_1^2(u)|_{\infty}\le c
\]
and
\[
\left| \frac{d}{d\tau}f_1^2\right|_{\infty}\le 
\frac{c}{|\lambda_m|}\|u-v\|
\]
at $\tau=0$, where $c$ is independent of $u$ on bounded sets in $M$.  
Use that $F_m(u)$ is an analytic functional of $u$ in $W$ to establish that the
curve  $u(\tau)$, $0\le \tau\le1$,
 remains in $M$ except for a countable number
of values of $\tau$. Begin with 
\[
f_{1,x}-uf_2=-i\lambda f_1
\]
and 
\[
f_{2,x}-\overline{u}f_1=i\lambda f_2.
\]
Multiply the equation for $f_1$ by $f_1$ and rewrite as
\[
\frac{1}{2}\partial_x f_1^2=uf_1f_2-i\lambda f_1^2.
\]
A similar calculation gives
\[
\frac{1}{2}\partial_x f_2^2-\overline{u}f_1f_2=i\lambda f_2^2.
\]
Next multiply the equation for $f_1$ by $f_2$ and combine with the equation for
$f_2$ multiplied by $f_1$ and find an expression for $\partial_x(f_1f_2)$. The
function 
\[
f_1f_2=\int_0^x \left(uf_2^2+\overline{u}f_1^2\right) dx
\]
satisfies the differential expression for $f_1f_2$. Substitute the above
identity for $f_1f_2$ into the preceding equation for $f_1^2$ and $f_2^2$ and
find
\begin{equation}
-\frac{1}{2}\partial_x f_1^2+u\int_0^x\left(uf_2^2+\overline{u}f_1^2\right)
dx=i\lambda f_1^2
\end{equation}
and
\[
\frac{1}{2}\partial_x f_2^2-\overline{u}\int_0^x
\left(uf_2^2+\overline{u}f_1^2\right) dx=i\lambda
f_2^2.
\]

Multiply (7) by $\exp(2i\lambda x)$ and rewrite the first equation
as 
\[
\partial\left(e^{2i\lambda x} f_1^2\right)=ue^{2i\lambda
x}\int^x\left(uf_2^2+\overline{u}f_1^2\right) dx.
\]
 Take the absolute value of this expression,
use the inequality
\[
\left| \partial\left(e^{2i\lambda x}f_1^2(x)\right)\right|\ge
\partial\left|e^{2i\lambda x}f_1^2(x)\right|
\]
 and then integrate in $x$ and find that
\[
\left| f_1^2(x)\right| \le \int^x\left| u\int^y\left(uf_2^2+
\overline{u}f_1^2\right)dy\right| dx\le|u|_{\infty}^2
\int^x\left(\left|f_1^2\right|+\left| f_2^2\right|\right)dx.
\]
 A similar calculation gives
\[
\left|f_2^2(x)\right|
\le |u|_{\infty}^2\int^x\left(\left|f_1^2\right|+\left|f_2^2\right|\right)dx.
\]
Combine and form 
$\left|f_1^2(x)\right|+\left| f_2^2(x)\right|$, use Gronwall and find that
\[
\left| f_1^2(x)\right| \le c,
\]
 where $c$ is independent of $u$ on bounded sets in $W$.
 
For $u(\tau)=\tau v+(1-\tau)u$, differentiate  (7) with respect to $\tau$ at
$\tau =0$ and find that 
\begin{gather*}
{\dot f}_1^2+\frac{1}{2i\lambda}\partial {\dot
f}_1^2=\frac{\dot\lambda}{\lambda}f_1^2-
\frac{\dot u}{i\lambda}\int^x\left(uf_2^2+\overline{u}f_1^2\right)dy\\
\phantom{{\dot f}_1^2+\frac{1}{2i\lambda}\partial {\dot f}_1^2=}{}
-\frac{u}{i\lambda}\int^x
\left(\dot u f_2^2+\dot{\overline{u}}f_1^2\right)dy-
\frac{u}{i\lambda}\int^x\left(u{\dot f}_2^2+\overline{u}{\dot
f}_1^2\right)dy
\end{gather*}
 and a similar expression for $\dot{f}_2^2$. Multiply the previous
identity by $\exp(x2i\lambda)$
and write the left hand side of the resulting expression as
\[
 \partial\left(e^{x2i\lambda}\dot{f}_1^2\right)/2i\lambda.
\]
 Next take the absolute value
of the expression, use the inequality
\[
c\partial \left|\dot{f}_1^2(x)\right|\le \partial
\left| e^{x2i\lambda}\dot{f}_1^2(x)/2i\lambda\right|
\le\left| \partial\left(e^{x2i\lambda}
\dot{f}_1^2(x)\right)/2i\lambda\right|,
\]
 where we have used that $\lambda$ is an analytic
functional of $u$ in $W$ and $c$ is an absolute constant. Integrate the resulting
expression in $x$ and use the estimates   
\begin{gather*}
\left| \int^xe^{y/2i\lambda}
\frac{{\dot u}}{i\lambda}\int^y\left(uf_2^2+\overline{u}f_1^2\right)\right|
 \le\int^x \frac{|\dot u|}{|\lambda|} 
\int^y\left|uf_2^2+\overline{u}f_1^2\right| \\
\qquad {}\le \frac{\|u-v\|}{|\lambda|}\|u\|\left(|f_2|_{\infty}^2+|f_1|_{\infty}^2\right) \le
\frac{c}{|\lambda|}\|u-v\|,\\
\left| \int^xe^{y/2i\lambda}\frac{u}{i\lambda}
\int^y\left(u\dot{f}_2^2+\overline{u}\dot{f}_1^2\right)\right| \le\int^x\frac{|u|}{|\lambda|} \int^y
\left| u\dot{f}_2^2+\overline{u}\dot{f}_1^2\right|\\
\qquad {}\le\frac{|u|_{\infty}^2}{|\lambda|}\int^x
\left(\left|\dot{f}_2^2\right|+\left|\dot{f}_1^2\right|\right)dy
 \le \frac{c}{|\lambda|}\int^x\left(\left|\dot{f}_2^2\right|+\left|\dot{f}_1^2\right| \right)dy,\\
\left|\frac{1}{\lambda}\int^x\left(if_1^2,u-v\right)f_1^2\right|\le
\frac{\|u-v\| |f_1|^2_{\infty}}{|\lambda|}\le\frac{c\|u-v\|}{|\lambda|},
\end{gather*} 
where we have used previous bounds and find
that
\[
\left| \dot{f}_1^2(x)\right|
\le \frac{c\|u-v\|}{|\lambda|}+
\frac{c}{|\lambda|}
\int^x \left(\left|\dot{f}_2^2\right|+\left|\dot{f}_1^2\right|\right)dy.
\] 
A similar
calculation gives  
\[
\left|\dot{f}_2^2(x)\right|
\le \frac{c\|u-v\|}{|\lambda|}+
\frac{c}{|\lambda|}\int^x\left(\left|\dot{f}_2^2\right|+\left|\dot{f}_1^2\right|\right)dy.
\]
 Combine  and form $\left|\dot{f}_1^2(x)\right|+\left|\dot{f}_2^2(x)\right|$
and use
Gronwall to obtain 
\[
\left|\dot{f}_1^2(x)\right|+\left|\dot{f}_2^2(x)\right|
\le \frac{c\|u-v\|}{|\lambda|},
\]
 where $c$ is independent of $u$ on bounded
sets of $M$.

 \medskip

\underline{\it Item 4.} For $u$ in $M$, the
normal vectors $G_{\lambda_m}(u)=if_{1,m}^2$ and the tangent vectors
$K_{\lambda_m}(u)= - f_{2,m}^2$ is a basis for $N_u$ and $T_u$ respectively and
$N_u\oplus T_u=L_2$. Each $G_u$ in $N_u$ is uniquely
$G_u=\sum\limits_{m=1}^{\infty}t_mG_{\lambda_m}(u)$ and  
\[
c_1|t|_{l_2}\le
\left\| \sum_{m=1}^{\infty}t_mG_{\lambda_m}(u)\right\|
\le c_2|t|_{l_2},
\]
 where $c_1$,
$c_2$ are independent of $u$ on bounded sets in  $M$. Each element $K_u$ in
the tangent space is uniquely represented as
$K_u=\sum\limits_{m=1}^{\infty}t_mK_{\lambda_m}(u)$ 
and  
\[
c_1|t|_{l_2}\le
\left\|\sum_{m=1}^{\infty}t_mK_{\lambda_m}(u)\right|
\le c_2|t|_{l_2},
\]
 where $c_1$,
$c_2$ are independent of $u$ on bounded sets in $M$.

We modify an idea of Borg [10] and 
establish this result by comparing the sequence $G_{\lambda_m}$
and $K_{\lambda_m}$ at $u$ with
the sequence at $u=0$. The work of McKean and Trubowitz~[11] 
used a~similar comparison in their
study of the basis properties of the normal and tangent space for the isospectral set of the
periodic Korteweg-de Vries equation.  For
$u=0$, we begin with the periodic and antiperiodic spectrum
$\lambda_{\pm}=n\pi$ and eigenfunctions
$f_{2m}=e^{ix\lambda_m}$. Then $G_{\lambda_m}=ie^{-i2x\lambda_m}$ and
$K_{\lambda_m}=-e^{-i2x\lambda_m}$. For
$u=0$, the closure in $L_2$ of the linear span of $G_{\lambda_m}$ and
$K_{\lambda_m}$ equals the closure in $L_2$ of the linear span of
\[
(a_n\sin(2n\pi x)+b_n\cos(2n\pi x))+i(a_n\cos(2n\pi x)+b_n\sin((2n\pi x)),
\]
 a basis for the space of complex valued functions that are square integrable.

For $u$=0, $G_{1,m}(0)=g_{1,m}^2$ is orthogonal and admits the estimate 
\[
cÕ_1|c|_{l_2}\le \left\|\sum_{m=1}^{\infty}c_mg_{1,m}^2\right\|\le
cÕ_2|c|_{l_2},
\]
 where $cÕ_1$ and $cÕ_2$ are absolute constants. We next
establish an apriori estimate of $f_{1,m}^2(u)$.

 For $u$ in $M$, define $T$
by 
\[
T\left(\sum_m c_m g_{1,m}^2\right)=\sum_mc_m
\left(f_{1,m}^2(u)-g_{1,m}^2\right).
\]
 Use of the
estimate in item 3 confirms that 
\[
\sum_m\left\| f_{1,m}^2(u)-g_{1,m}^2\right\|^2<\infty
\]
and establishes that $T$ is a bound linear operator and is Hilbert Schmidt. If
$f_{1,m}^2(u)$ is minimal then $(I+T)$
 is invertible and 
\[
\left| (I+T)^{-1}\sum_mc_mf_{1,m}^2(u_1)\right|
=\sum_m|c_m|^2\left\| g_{1,m}^2\right\|^2\le
c' |c|_{l_2},
\]
 where we have used that $\left(g_{1,m}^2,g_{1,n}^2\right)=0$,
$m\ne n$, the
estimate of item 3, and $c'$ is independent of $u$ on bounded sets from $M$. It
follows that 
\begin{equation}
\left\| \sum_mc_mf_{1,m}^2(u)\right\|
\ge c_1|c|_{l_2},
\end{equation}
where $c_1$ is independent of $u$ on bounded sets in $M$.

 We use the apriori
estimate and establish that the basis $g_{1,m}^2$ may be continuously extended to
a basis $f_{1,m}^2$. Let $u(\tau)=\tau u$ and use that $F_m(u)$ is an
analytic functional of~$u$ in $W$ to establish that the curve $u(\tau)$,
$0\le\tau\le 1$ remains in $M$ except for a countable number of $\tau$ and write
\[
f_{1,m}^2(\tau u)-g_{1,m}^2=\int_0^{\tau}\frac{d}{d\tau} f_{1,m}^2(su)\,ds.
\]
 Use the estimate of item 3 and find 
\[
\left\| f_{1,m}^2(\tau u)-g_{1,m}^2\right\|
\le\tau\left\|\frac{d}{ds}f_{1,m}^2(su)\right\|
\le\frac{\tau}{\|\lambda_m\|}\, c\|u\|,
\]
where $c$ is independent of $u$ on bounded sets of $M$. Use (7), the previous
bound, and select $\tau=\tau_1$ independently of $u$ on bounded sets of $M$ and
find $u_1=\tau_1 u$ in $M$, where 
\[
\left\| \sum_mc_mf_{1,m}^2(u_1)\right\|
\ge\frac{c_1}{2}|c|_{l_2}.
\]
 This estimate confirms that $f_{1,m}^2(u_1)$ is a minimal
sequence and that the closure in $L_2$ of the linear span of $f_{1,m}^2(u_1)$
equals the closure in $L_2$ of the linear span of $g_{1,m}^2$. Use of the
apriori bound (8) gives the estimate
\[
\left\| \sum_mc_mf_{1,m}^2(u_1)\right\|
\ge c_1 |c|_{l_2}.
\]
 Iteration of this
construction gives the result. 

\medskip

\underline{\it Item 5.} For $u$ in $M$, $dG_u$ admits the estimate
$(dG_uK_u,K_uÕ)/\|K_u\|\|K_uÕ\|\le c\|G_u\|$, where $c$ is independent of
$u$ on bounded sets in $M$. For the direction $v$ transverse to $M_u$ at~$u$,
$dG_m(u)$ satisfies the estimate $\|dG_m(u)v\|\le c_m\|v\|$,
 where $c_m$ is
square summable uniformly in $u$ and $v$ on bounded sets in $M$. 

 Use of (3) and a direct
calculation shows that $\lambda_m / m \rightarrow 1$ uniformly in $u$ on
bounded sets in $M$ as $m\rightarrow \infty$. Combine this result and the
estimates of item 3 to obtain the estimates of this section.

\medskip

This completes the proof of the Theorem.

\label{schwarz-lastpage}
\end{document}